\documentclass[12pt,a4paper]{article}
\usepackage[utf8]{inputenc}
\usepackage{amsmath}
\usepackage{amsfonts}
\usepackage{amssymb}
\usepackage{graphicx}
\usepackage[left=4cm,right=4cm,top=3cm,bottom=3cm]{geometry}
\usepackage{amsthm}
\newtheorem{thm}{Theorem}[section]

\newtheorem{defn}[thm]{Definition}

\newtheorem{exm}[thm]{Example}
\usepackage[colorlinks,urlcolor=blue]{hyperref}
\usepackage{enumerate}
\usepackage{multirow}
\numberwithin{equation}{subsection}

\title{Ranking of Intuitionistic Fuzzy Numbers by New Distance Measure}

\author{Debaroti Das$^{1}$,P.K.De$^{2*}$\\ $^{1,2}$Department of Mathematics\\ NIT Silchar,788010, Assam, India\\ 
Email: deboritadas1988@gmail.com$^{1}$, pijusde@gmail.com$^{2}$}
\date{}

\begin{document}
\maketitle

\begin{abstract}
Ranking of intuitionsitic fuzzy number plays a vital role in decision making and other intuitionistic fuzzy applications. In this paper, we propose a new ranking method of intuitionistic fuzzy number based on distance measure. We first define a distance measure for interval numbers based on $L_{p}$ metric and further generalize the idea for intuitionistic fuzzy number by forming interval with their respective value and ambiguity indices. Finally, some comparative results are given in tabular form.

\end{abstract}

\textbf{keywords:}Trapezoidal intuitionistic fuzzy number(TRIFN), Triangular intuitionistic fuzzy number(TIFN), value and ambiguity index, ranking , distance measure




\section{Introduction}
\label{}
Zadeh \cite{za} introduced Fuzzy set theory in 1965. Later on Atanassov generalized the concept of fuzzy set and introduced the idea of Intuitionistic fuzzy set \cite{at}. Intuitionistic fuzzy sets (IFSs) are characterized by three functions each of them expressing degrees of membership, non-membership, and indeterminacy \cite{atg}. The notions of intuitionistic fuzzy numbers in different context were studied in [\cite{lm}- \cite{lra},\cite{mir},\cite{wan},\cite{wz},\cite{wu}] and applied in multi criteria decision making problems [\cite{lm}-\cite{ye}]. The ranking of intuitionistic fuzzy numbers plays a vital role in intuitionistic fuzzy decision-making. In last few years, different methods for ranking IFNs have been introduced[\cite{hua}-\cite{lm},\cite{lra},\cite{mir},\cite{xuc}]. Szmidt and kacpryzk introduced the hamming distance between intuitionistic fuzzy sets and proposed a similarity measure based on the distance\cite{kp}. In \cite{skp},Szmidt and kacpryzk suggested some methods for measuring distances between intuitionistic fuzzy sets that are well known generalization of Hammimg distance, Euclidean distance and their normalised counterparts. Few methods for measuring distances between IFS/IVIFS based on Hausdorff metric were suggested by P.Grzegorzewski\cite{grz}. Hung and Yang\cite{hung} also introduced a different method for measuring distances and ranking of intuitionistic fuzzy number based on Hausdorff distance. Huang et al.\cite{hua} proposed the distance measures for IFSs that are generalization of Minkowski's distance of fuzzy number. Based on Hausdorff distance, Wan \cite{wn} defined Hamming distance and Euclidean distance for trapezoidal intuitionistic fuzzy numbers. Mitchell\cite{mir} introduced ranking of IFNs. Nayagam et al.\cite{lm} described IFNs of special type and introduced a method of IF scoring of IFNs as a generalization of of Chen and Hwang's scoring for ranking IFNs\cite{chen}. Wang and Zhang\cite{wz} defined the trapezoidal Intuitionistic fuzzy number(TRIFN) and gave ranking method which transformed the ranking of TRIFNs in to ranking of interval numbers. Xu\cite{xuzs} and Xu and Chen \cite{xuc} proposed score function and accuracy function to rank intervalvalued intuitionistic fuzzy numbers. Ye\cite{ye} also proposed a novel accuracy function to rank interval-valued intuitionistic fuzzy numbers.  Li \cite{lra} defined the value index and the ambiguity index of triangular intuitionistic fuzzy numbers and developed a ratio ranking method for solving multi-attribute decision-making problems. Wan and Li \cite{wan} focused on multi-attribute decision making (MADM) problems in which the attribute values are expressed with TIFNs and the information on attribute weights is incomplete, which are solved by developing a new decision method based on possibility mean and variance of TIFNs. Wang\cite{wang} gave the definition of intuitionistic trapezoidal fuzzy number. Wang et al.\cite{wli} also developed a methodology for solving multiattribute decision making (MADM) with both ratings of alternatives on attributes and weights being expressed with IVIF sets. Further, Wang and Zhang\cite{wz} developed the hamming distance for intuitionistic trapezoidal fuzzy numbers. Recently,for multi-criteria group decision making problems with intuitionistic linguistic information, Wang et al.\cite{wwy}, defined a new score function and a new accuracy function of intuitionistic linguistic numbers, and proposed a simple approach for the comparison between two intuitionistic linguistic numbers.  In the present article,authors suggested a new method of ranking of intuitionistic fuzzy number which is based on $L_{p}$ metric.
\\\\\hspace*{.4in} The remainder of the paper is organised as follows.Section 2,contains preliminaries followed by Section 3 which is consists of distance measure for interval numbers and distance measure for intuitionistic fuzzy numbers. In Section 4, we present the proposed ranking method. The last two sections,Section 5 and section 6 include few illustrative examples to compare the proposed method with other methods and conclusive remarks respectively.

\section{Preliminaries}
We shall start this section with some basic notations and definition related  to intuitionistic fuzzy numbers. 

\begin{defn}
\cite{wz}
 A TRIFN $\tilde{a}=\langle (a_{1},a_{2},a_{3},a_{4});w_{\tilde{a}},u_{\tilde{a}} \rangle$ is a special Intuitionistic Fuzzy set on a set of real number $\mathbb{R}$, whose membership function and non membership function are defined as follows:
\[\mu_{\tilde{a}}(x)= \left\{\begin{array}{cc}
\dfrac{(x-a_{1})}{(a_{2}-a_{1})} w_{\tilde{a}} & a_{1}\leq x< a_{2} \\ 
w_{\tilde{a}} &  a_{2}\leq x\leq a_{3}\\
\dfrac{(a_{4}-x)}{(a_{4}-a_{3})} w_{\tilde{a}}& a_{3}< x\leq a_{4} \\ 
0 & a_{4}<x, a_{1}>x
\end{array}
\right. 
\]\\ and
\[\nu_{\tilde{a}}(x)= \left\{\begin{array}{cc}
\dfrac{(a_{2}-x)+u_{\tilde{a}}(x-a_{1})}{(a_{2}-a_{1})} & a_{1}\leq x<a_{2} \\ 
u_{\tilde{a}} &  a_{2}\leq x\leq a_{3} \\ 
\dfrac{(x-a_{3})+u_{\tilde{a}}(a_{4}-x)}{(a_{4}-a_{3})} & a_{3}< x\leq a_{4} \\ 
1 & a_{4}<x ,a_{1}>x
\end{array}
\right. 
\]
respectively.
\begin{figure}[hbtp]
\centering
\includegraphics[scale=1]{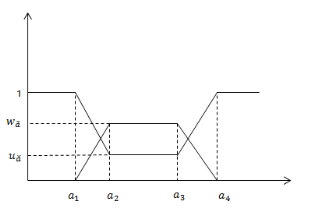}
\caption{Trapezoidal Intuitionistic fuzzy numbers}
\end{figure}

 The values $w_{\tilde{a}}$ and $u_{\tilde{a}}$ represents the maximum degree of membership and minimum degree of non membership, respectively, such that the conditions $0 \leq w_{\tilde{a}}\leq 1$ , $0\leq u_{\tilde{a}}\leq 1$ and $0 \leq w_{\tilde{a}}+u_{\tilde{a}} \leq 1$ are satisfied. The parameters $w_{\tilde{a}}$ and $u_{\tilde{a}}$ reflects the confidence level and non confidence level of the TRIFN $\tilde{a}=\langle (a_{1},a_{2},a_{3},a_{4});w_{\tilde{a}},u_{\tilde{a}} \rangle$,respectively.
\\\hspace*{.5in} Let  $ \pi_{\tilde{a}}(x)=1-\mu_{\tilde{a}}(x)-\nu_{\tilde{a}}(x) $ which is called an IF index of an element x in $\tilde{a}$.It is the degree of the indeterminacy membership of the element x in $\tilde{a}$.
\end{defn}

\begin{defn}\cite{lra}
Let $\tilde{a}=\langle (a_{1},a_{2},a_{3},a_{4});w_{\tilde{a}},u_{\tilde{a}} \rangle$ and $\tilde{b}=\langle (b_{1},b_{2},b_{3},b_{4});w_{\tilde{b}},u_{\tilde{b}} \rangle$ be two TRIFNs and $\lambda$ be a real number. The arithmetical operations are listed as follows:
\\\\\small{$\tilde{a}+\tilde{b}=\langle(a_{1}+b_{1},a_{2}+b_{2},a_{3}+b_{3},a_{4}+b_{4});min\lbrace w_{\tilde{a}},w_{\tilde{b}}\rbrace,max\lbrace u_{\tilde{a}},u_{\tilde{b}}\rbrace\rangle$}
\\\small{ $\tilde{a}-\tilde{b}=\langle(a_{1}-b_{4},a_{2}-b_{3},a_{3}-b_{2},a_{4}+b_{1});min\lbrace w_{\tilde{a}},w_{\tilde{b}}\rbrace,max\lbrace u_{\tilde{a}},u_{\tilde{b}}\rbrace\rangle$} 
\\\small{ $\\\tilde{a}\tilde{b} =\left\{ \begin{array}{lll}
\langle(a_{1}b_{1},a_{2}b_{2},a_{3}b_{3},a_{4}b_{4});min\lbrace w_{\tilde{a}},w_{\tilde{b}}\rbrace,max\lbrace u_{\tilde{a}},u_{\tilde{b}}\rbrace\rangle & if \tilde{a}>0 & and \tilde{b}>0 \\ 
\langle(a_{1}b_{4},a_{2}b_{3},a_{3}b_{2},a_{4}b_{1});min\lbrace w_{\tilde{a}},w_{\tilde{b}}\rbrace,max\lbrace u_{\tilde{a}},u_{\tilde{b}}\rbrace\rangle & if \tilde{a}<0 & and \tilde{b}>0 \\ 
\langle(a_{4}b_{4},a_{3}b_{3},a_{2}b_{2},a_{1}b_{1});min\lbrace w_{\tilde{a}},w_{\tilde{b}}\rbrace,max\lbrace u_{\tilde{a}},u_{\tilde{b}}\rbrace\rangle & if \tilde{a}<0 & and \tilde{b}<0.
\end{array} 
\right.$ }
\\\small{$\lambda \tilde{a}=\left\{ \begin{array}{ll}
\langle(\lambda a_{1},\lambda a_{2},\lambda a_{3},\lambda a_{4});w_{\tilde{a}},u_{\tilde{a}} \rangle & if \lambda>0\\ 
\langle(\lambda a_{4},\lambda a_{3},\lambda a_{2},\lambda a_{1});w_{\tilde{a}},u_{\tilde{a}} \rangle & if \lambda<0
\end{array}
\right. $}
\\\small{$\tilde{a}^{-1}= \langle(1\diagup a_{1},1\diagup a_{2}, 1\diagup a_{3}, 1\diagup a_{4}),w_{\tilde{a}},u_{\tilde{a}}\rangle$}
\end{defn}

\textbf{$ \alpha$-cut sets and $\beta$-cut sets of TRIFNs:}

\begin{defn}\cite{lra}
A $\alpha$-cut set of a TRIFN $\tilde{a}=\langle (a_{1},a_{2},a_{3},a_{4});w_{\tilde{a}},u_{\tilde{a}} \rangle$ is a crisp sub set of $\mathbb{R}$, defined as $\tilde{a}_{\alpha}=\lbrace x\vert \mu_{\tilde{a}}(x)\geq \alpha \rbrace$ where $0 \leq \alpha \leq w_{\tilde{a}}$.
\end{defn}

\begin{defn}\cite{lra}
 A $\beta$-cut set of a TRIFN $\tilde{a}=\langle (a_{1},a_{2},a_{3},a_{4});w_{\tilde{a}},u_{\tilde{a}} \rangle$ is a crisp sub set of $\mathbb{R}$, defined as $\tilde{a}_{\beta}=\lbrace x\vert \nu_{\tilde{a}}(x)\leq \beta \rbrace$ where $ u_{\tilde{a}} \leq \beta \leq 1$.
\end{defn}
It can be easily followed from the definition of TRIFN and Definition 2.3 and Definition 2.4 that $\tilde{a}_{\alpha}$ and $\tilde{a}_{\beta}$ are both closed sets and are denoted by $\tilde{a}_{\alpha}=[L_{\tilde{a}}(\alpha),R_{\tilde{a}}(\alpha)]$ and $\tilde{a}_{\beta}=[L_{\tilde{a}}(\beta),R_{\tilde{a}}(\beta)]$ respectively. The respective values of $\tilde{a}_{\alpha}$ and $\tilde{a}_{\beta}$ are calculated as follows:
\\$[L_{\tilde{a}}(\alpha),R_{\tilde{a}}(\alpha)]=[a_{1}+\dfrac{\alpha(a_{2}-a_{1})}{w_{\tilde{a}}},a_{4}-\dfrac{\alpha(a_{4}-a_{3})}{w_{\tilde{a}}}]$ and\\ $[L_{\tilde{a}}(\beta),R_{\tilde{a}}(\beta)]=[\dfrac{(1-\beta)a_{2}+(\beta-u_{\tilde{a}})a_{1}}{1-u_{\tilde{a}}},\dfrac{(1-\beta)a_{3}+(\beta-u_{\tilde{a}})a_{4}}{1-u_{\tilde{a}}}]$

\subsection{Definition of Value and Ambiguity Index for TRIFNs}
\label{sec:2.1}
 Value and Ambiguity of a triangular fuzzy number were generalised by D.F.Li \cite{lra}. The value and ambiguity of a TRIFN have been defined by the author \cite{pk} in analogue with value and ambiguity of TIFN defined by D.F.Li\cite{lra}.
 
 \begin{defn} \cite{pk} Let $\tilde{a}_{\alpha}$ and $\tilde{a}_{\beta}$ be an $\alpha$-cut set and a $\beta$-cut set of a TRIFN $\tilde{a}=\langle (a_{1},a_{2},a_{3},a_{4});w_{\tilde{a}},u_{\tilde{a}} \rangle$, respectively. Then the values of the membership function $\mu_{\tilde{a}}(x)$ and the non-membership function $\nu_{\tilde{a}}(x)$ for the TRIFN $\tilde{a}$ are defined as follows:
\[V_{\mu}(\tilde{a})=\int_{0}^{w_{\tilde{a}}} \dfrac{L_{\tilde{a}}(\alpha)+R_{\tilde{a}}(\alpha)}{2}f(\alpha)d\alpha\]
\[V_{\nu}(\tilde{a})=\int_{u_{\tilde{a}}}^{1} \dfrac{L_{\tilde{a}}(\beta)+R_{\tilde{a}}(\beta)}{2}g(\beta)d\beta \]
respectively,where the function $f(\alpha)$ is a non-negative and non-decreasing function on the interval $[0,w_{\tilde{a}}]$ with $f(0)=0$ and $\int_{0}^{w_{\tilde{a}}}f(\alpha)d\alpha= w_{\tilde{a}}$; the function $g(\beta)$ is a non-negative and non-increasing function on the interval $[u_{\tilde{a}},1]$ with $g(1)=0$ and $\int_{w_{\tilde{a}}}^{1}g(\beta)d\beta= 1-u_{\tilde{a}}$.
\end{defn}

\begin{defn} \cite{pk} Let $\tilde{a}_{\alpha}$ and $\tilde{a}_{\beta}$ be an $\alpha$-cut set and a $\beta$-cut set of a TRIFN $\tilde{a}=\langle (a_{1},a_{2},a_{3},a_{4});w_{\tilde{a}},u_{\tilde{a}} \rangle$, respectively. Then the ambiguities of the membership function $\mu_{\tilde{a}}(x)$ and the non-membership function $\nu_{\tilde{a}}(x)$ for the TRIFN $\tilde{a}$ are defined as follows:

\[A_{\mu}(\tilde{a})=\int_{0}^{w_{\tilde{a}}} (R_{\tilde{a}}(\alpha)-L_{\tilde{a}}(\alpha))f(\alpha)d\alpha\]
\[A_{\nu}(\tilde{a})=\int_{u_{\tilde{a}}}^{1} (R_{\tilde{a}}(\beta)-L_{\tilde{a}}(\beta))g(\beta)d\beta \]
respectively. 
\end{defn}
For rest of the paper we shall choose $f(\alpha)=\dfrac{2\alpha}{w_{\tilde{a}}}$,$\alpha \in [0,w_{\tilde{a}}]$ and $g(\beta)=\dfrac{2(1-\beta)}{1-u_{\tilde{a}}}$, $\beta \in [u_{\tilde{a}},1]$ and the value of membership and non-membership of a TRIFN $\tilde{a}$ have been derived \cite{pk} as follows : 
\\The value of the membership function of a TRIFN $\tilde{a}$ is calculated as follows:

\begin{align*}
V_{\mu}(\tilde{a}) &=\int_{0}^{w_{\tilde{a}}}\left[ a_{1}+\dfrac{\alpha(a_{2}-a_{1})}{w_{\tilde{a}}}+a_{4}-\dfrac{\alpha(a_{4}-a_{3})}{w_{\tilde{a}}}\right]\dfrac{\alpha}{w_{\tilde{a}}}d\alpha \\
 &=  \left[\dfrac{a_{1}+a_{4}}{2w_{\tilde{a}}}\alpha^{2}\right]_{0}^{w_{\tilde{a}}}+\left[\dfrac{(a_{2}-a_{1}-a_{4}+a_{3})}{3(w_{\tilde{a}})^{2}}\alpha^{3}\right]_{0}^{w_{\tilde{a}}} \\
 &=  \dfrac{a_{1}+a_{4}+2(a_{2}+a_{3})}{6} w_{\tilde{a}}
\end{align*}

In a similar way value of non-membership can be evaluated as follows:
\begin{align*}
V_{\nu}(\tilde{a}) &=\int_{u_{\tilde{a}}}^{1}\left[\dfrac{(1-\beta)a_{2}+(\beta-u_{\tilde{a}})a_{1}}{1-u_{\tilde{a}}}+\dfrac{(1-\beta)a_{3}+(\beta-u_{\tilde{a}})a_{4}}{1-u_{\tilde{a}}} \right] \dfrac{1-\beta}{1-u_{\tilde{a}}}d\beta \\
&=\int_{u_{\tilde{a}}}^{1}\dfrac{(a_{2}+a_{3}-a_{1}-a_{4})(1-\beta)^{2}+ (1-u_{\tilde{a}})(a_{1}+a_{4})(1-\beta)}{(1-u_{\tilde{a}})} d\beta\\ 
&= -\left[\dfrac{(a_{2}+a_{3}-a_{1}-a_{4})(1-\beta)^{3}}{3(1-u_{\tilde{a}})^{2}}\right]_{u_{\tilde{a}}}^{1}-\left[\dfrac{(a_{1}+a_{4})(1-u_{\tilde{a}})(1-\beta)^{2}}{2(1-u_{\tilde{a}})^{2}} \right]_{u_{\tilde{a}}}^{1} \\
&= \dfrac{[a_{1}+a_{4}+2(a_{2}+a_{3})](1-u_{\tilde{a}})}{6}
\end{align*}

 With the condition that $0 \leq w_{\tilde{a}}+ u_{\tilde{a}} \leq 1$ ,it follows that $V_{\mu}(\tilde{a}) \leq V_{\nu}(\tilde{a})$. Thus, the values of membership and non-membership functions of a TRIFN $\tilde{a}$ may be concisely expressed as an interval $[V_{\mu}(\tilde{a}), V_{\nu}(\tilde{a})]$.
 
The ambiguity of the membership function of a TRIFN $\tilde{a}$ is calculated as follows \cite{pk}:
\begin{align*}
A_{\mu}(\tilde{a})&=\int_{0}^{w_{\tilde{a}}}\left[ a_{4}-\dfrac{\alpha(a_{4}-a_{3})}{w_{\tilde{a}}}-a_{1}-\dfrac{\alpha(a_{2}-a_{1})}{w_{\tilde{a}}}\right]\dfrac{2\alpha}{w_{\tilde{a}}}d\alpha\\
&= \left[\dfrac{a_{4}-a_{1}}{w_{\tilde{a}}}\alpha^{2}\right]_{0}^{w_{\tilde{a}}}-2\left[\dfrac{(a_{2}-a_{1}+a_{4}-a_{3})}{3(w_{\tilde{a}})^{2}}\alpha^{3}\right]_{0}^{w_{\tilde{a}}}\\
&=\dfrac{(a_{4}-a_{1})-2(a_{2}-a_{3})}{3}w_{\tilde{a}}
\end{align*}

Similarly, the ambiguity of the non-membership function of a TRIFN $\tilde{a}$ is given by:

\begin{align*}
A_{\nu}(\tilde{a})&=\int_{u_{\tilde{a}}}^{1}\left[\dfrac{(1-\beta)a_{3}+(\beta-u_{\tilde{a}})a_{4}}{1-u_{\tilde{a}}}-\dfrac{(1-\beta)a_{2}+(\beta-u_{\tilde{a}})a_{1}}{1-u_{\tilde{a}}} \right] \dfrac{2(1-\beta)}{1-u_{\tilde{a}}}d\beta \\
&=\int_{u_{\tilde{a}}}^{1}\dfrac{2\left[-(a_{2}-a_{3}-a_{1}+a_{4})(1-\beta)^{2}+ (1-u_{\tilde{a}})(a_{4}-a_{1})(1-\beta)\right]}{(1-u_{\tilde{a}})^2} d\beta \\
&=\left[\dfrac{2(a_{4}-a_{1}+a_{2}-a_{3})(1-\beta)^{3}}{3(1-u_{\tilde{a}})^{2}}\right]_{u_{\tilde{a}}}^{1}- \left[\dfrac{(a_{4}-a_{1})(1-u_{\tilde{a}})(1-\beta)^{2}}{(1-u_{\tilde{a}})^{2}}\right]_{u_{\tilde{a}}}^{1} \\
&=\dfrac{(a_{4}-a_{1})-2(a_{2}-a_{3})}{3}(1-u_{\tilde{a}})
\end{align*}
 With the condition that $0 \leq w_{\tilde{a}}+ u_{\tilde{a}} \leq 1$ ,it follows that $A_{\mu}(\tilde{a}) \leq A_{\nu}(\tilde{a})$. Thus, the ambiguities of membership and non-membership functions of a TRIFN $\tilde{a}$ may be concisely expressed as an interval $[A_{\mu}(\tilde{a}), A_{\nu}(\tilde{a})]$.

\subsection{Value and Ambiguity Index for TRIFNs }
\label{sec 2.2}
\begin{defn} \cite{pk} Let $\tilde{a}=\langle (a_{1},a_{2},a_{3},a_{4});w_{\tilde{a}},u_{\tilde{a}} \rangle$ be a TRIFN. A value index and ambiguity index for the TRIFN $\tilde{a}$ are defined as follows:
\[V(\tilde{a},\lambda)= V_{\mu}(\tilde{a})+\lambda(V_{\nu}(\tilde{a})-V_{\mu}(\tilde{a}))
\hspace{.5in} A(\tilde{a},\lambda)= A_{\nu}(\tilde{a})-\lambda(A_{\nu}(\tilde{a})-A_{\mu}(\tilde{a}))\]
\\respectively, where $\lambda \in [0,1]$ is a weight which represents the decision maker's preference information. Limited to the above formulation,the choice $\lambda=\dfrac{1}{2}$ appears to be a reasonable one. One can choose $\lambda$ according to the suitability of the subject. $\lambda\in [0,\dfrac{1}{2})$ indicates decision maker's pessimistic attitude towards uncertainty while $\lambda\in (\dfrac{1}{2},1]$ indicates decision maker's optimistic attitude towards uncertainity.
\\\\With our choice $\lambda=\dfrac{1}{2}$, the value and ambiguity indices for TRIFN reduces to the following:
\[V(\tilde{a},\dfrac{1}{2})= \dfrac{V_{\mu}(\tilde{a})+V_{\nu}(\tilde{a})}{2}
\hspace{.5in} A(\tilde{a},\dfrac{1}{2})= \dfrac{A_{\mu}(\tilde{a})+A_{\nu}(\tilde{a})}{2}\]
\\Since we are taking $\lambda=\dfrac{1}{2}$ throughout the paper,instead of $V(\tilde{a},\dfrac{1}{2})$ and $A(\tilde{a},\dfrac{1}{2})$ we shall write $V(\tilde{a})$ and$A(\tilde{a})$ respectively.
\end{defn}

\section{Proposed distance measure}
\label{sec:3}
In this section, we define the distance measure for two interval numbers. Based of the newly defined interval number distance measure we shall define distance measure for trapezoidal intuitionistic fuzzy numbers.
\subsection{Distance measure between two interval numbers}
\label{sec:3.1}
\begin{defn}Let $[a,b]$ and $[c,d]$ be two intervals and let $f(x)=(a-b)x+b$ and $g(x)=(c-d)x+d$. The distance between two intervals is denoted and defined by 
\begin{align} 
\label{eq:1}
d^{(p)}([a,b],[c,d])= \|f(x)-g(x)\|_{L_{p}} 
\end{align}
where $\|.\|$ is the usual norm in the $L_{p}$ space on the interval $[0,1]$ $(p>1)$.
\end{defn}

The advantage of this distance measure is that it takes into account every point on the both intervals with convex combination of lower and upper bound values of the both intervals.
\\\hspace*{.3in} Much identical definition was followed by Allahviranloo and Firozja \cite{al}, but application of their proposed distance measure $d_{TMI}^{(p)}$ to intuitionistic fuzzy setting would make the scenario more complex. To make this distance measure more convinient to be applicable in case of intuitionistic fuzzy setting,the authors have defined it in a much more easier way.
\subsection{Metric properties}
\label{sec:3.2}
The distance measure as defined in \ref{eq:1} satisfies the following metric properties:
\begin{itemize}
\item[1.] $d^{(p)}([a,b],[c,d]) \geq 0$
\item[2.] $d^{(p)}([a,b],[c,d])=0 if and only if [a,b]=[c,d]$
\item[3.] $d^{(p)}([a,b],[c,d])=d^{(p)}([c,d],[a,b])$
\item[4.] $d^{(p)}([a,b],[c,d])+d^{(p)}([c,d],[e,z])\geq d^{(p)}([a,b],[e,z])$
\end{itemize}
\subsection{Various Properties of the distance measure $d^{(p)}$}
\label{sec:3.3}
Let $a,b,c,d, e,z$ are real numbers. Propositions given by Allahviranloo and Firozja \cite{al} can be reformulated as follows:
\small{
\begin{enumerate}
\item[Proposition 1:] $d^{(p)}([a+\lambda,b+\lambda],[c+\lambda,d+\lambda])=d^{(p)}([a,b],[c,d])$. 
\item[Proposition 2:] $d^{(p)}([\lambda a,\lambda b],[\lambda c,\lambda d])=|\lambda|d^{(p)}([a,b],[c,d])$.
\item[Proposition 3:] If $p=2$ then $d^{(p)}([a,b],[c,d])=\sqrt{\dfrac{1}{3}[(a-c)^{2}+(b-d)^{2}+(a-c)(b-d)]}$ .
\item[Proposition 4:] If $p=2$ then $d^{(p)}([a,b],c)=\sqrt{\dfrac{1}{3}[(a-c)^{2}+(b-c)^{2}+(a-c)(b-c)]}$ .
\item[Proposition 5:] If $p=2$ then $d^{(p)}([a,b],0)=\sqrt{\dfrac{1}{3}[a^{2}+b^{2}+ab]}$.
\item[Proposition 6:] If $p=3$ then $d^{(p)}([a,b],[c,d])= \left( \dfrac{1}{4}[(a-c)^{3}+(b-d)^{3}+(a-c)^{2}(b-d)+(a-c)(b-d)^{2}]\right)^{1/3} $.
\item[Proposition 7:] If $p=3$ then $d^{(p)}([a,b],0)= \left( \dfrac{1}{4}[a^{3}+b^{3}+a^{2}b +ab^{2}]\right)^{1/3} $.
\end{enumerate}}

\subsection{Distance measure for Intuitionistic fuzzy numbers}
\label{sec:3.4}
\begin{defn}
For any two trapezoidal intutionistic fuzzy numbers $\tilde{a}=\langle (a_{1},a_{2},a_{3},a_{4});w_{\tilde{a}},u_{\tilde{a}} \rangle$ and $\tilde{b}=\langle (b_{1},b_{2},b_{3},b_{4});w_{\tilde{b}},u_{\tilde{b}} \rangle$, we can evaluate value and ambiguity indices $V(\tilde{a}),V(\tilde{b}),A(\tilde{a}),A(\tilde{b})$ respectively. Based on the proposed distance measure, the distance measure for two TRIFNs have been defined as follows:

\begin{eqnarray*}
D(\tilde{a},\tilde{b})&=&
\begin{array}{cc}
d^{(p)}\left( [A(\tilde{a}), V(\tilde{a})],[A(\tilde{b}), V(\tilde{b})]\right) : V(\tilde{a}) \geq 0 \\
d^{(p)}\left( [V(\tilde{a}), A(\tilde{a})],[V(\tilde{b}), A(\tilde{b})]\right) : V(\tilde{a}) < 0
\end{array}\\
&=& \begin{array}{cc}
\left(\int_{0}^{1}{|A(\tilde{a})x+(1-x)V(\tilde{a})-A(\tilde{b})x-(1-x)V(\tilde{b})|}^{p}dx \right)^{1/p} :  V(\tilde{a}) \geq 0 \\
\left(\int_{0}^{1}{|V(\tilde{a})x+(1-x)A(\tilde{a})-V(\tilde{b})x-(1-x)A(\tilde{b})|}^{p}dx \right)^{1/p} :  V(\tilde{a}) < 0
\end{array}
\end{eqnarray*}

\end{defn}

\subsection{Metric properties of distance measure for Intuitionistic fuzzy numbers}
\label{sec:3.5}
The Proposed distance measure is a metric distance and satisfies the metric properties:
\begin{itemize}
\item[1.] $D(\tilde{a},\tilde{b}) \geq 0$.
\begin{proof}
 According to definition \\$D(\tilde{a},\tilde{b})=\left(\int_{0}^{1}{|A(\tilde{a})x+(1-x)V(\tilde{a})-A(\tilde{b})x-(1-x)V(\tilde{b})|}^{p}dx \right)^{1/p}$.
 \\The proof follows from the fact that\\ $|A(\tilde{a})x+(1-x)V(\tilde{a})-A(\tilde{b})x-(1-x)V(\tilde{b})| \geq 0$.
\end{proof}
\item[2.] $D(\tilde{a},\tilde{b})=0$ if and only if $\tilde{a}=\tilde{b}$.
\begin{proof}
 The proof follows from the fact that $\tilde{a}=\tilde{b}$ if and only if $\tilde{a}_{\alpha}=\tilde{b}_{\alpha}$ and $\tilde{a}_{\beta}=\tilde{b}_{\beta}$. Consequently, $\tilde{a}=\tilde{b}$ if and only if $V(\tilde{a})= V(\tilde{b})$ and $A(\tilde{a})=A(\tilde{b})$.
\end{proof}
\item[3.] $D(\tilde{a},\tilde{b})=D(\tilde{b},\tilde{a})$.
\begin{proof} 
\begin{align*}
D(\tilde{a},\tilde{b})&=\left(\int_{0}^{1}{|A(\tilde{a})x+(1-x)V(\tilde{a})-A(\tilde{b})x-(1-x)V(\tilde{b})|}^{p} dx \right)^{1/p}\\
&=\left(\int_{0}^{1}{|-(A(\tilde{b})x+(1-x)V(\tilde{b})-A(\tilde{a})x-(1-x)V(\tilde{a}))|}^{p} dx \right)^{1/p}\\
&=\left(\int_{0}^{1}{|-(A(\tilde{b})x+(1-x)V(\tilde{b})-A(\tilde{a})x-(1-x)V(\tilde{a}))|}^{p} dx \right)^{1/p}\\
&=\left(\int_{0}^{1}{|A(\tilde{b})x+(1-x)V(\tilde{b})-A(\tilde{a})x-(1-x)V(\tilde{a})|}^{p}dx \right)^{1/p}\\
&=D(\tilde{b},\tilde{a})
\end{align*} 
\end{proof}
\item[4.] $D(\tilde{a},\tilde{c})+D(\tilde{c},\tilde{b})\geq D(\tilde{a},\tilde{b})$
\begin{proof} $A(.)x+(1-x)V(.) \in L_{p}(0,1)$ for all TRIFNs.

\begin{align*}
D(\tilde{a},\tilde{c})+D(\tilde{c},\tilde{b})&=\left(\int_{0}^{1}{|A(\tilde{a})x+(1-x)V(\tilde{a})-A(\tilde{c})x-(1-x)V(\tilde{c})|}^{p}dx \right)^{1/p} 
\nonumber \\
 &\qquad {}
+ \left(\int_{0}^{1}{|A(\tilde{c})x+(1-x)V(\tilde{c})-A(\tilde{b})x-(1-x)V(\tilde{b})|}^{p}dx \right)^{1/p}\\
&\geq \left(\int_{0}^{1} |A(\tilde{a})x+(1-x)V(\tilde{a})-A(\tilde{c})x-(1-x)V(\tilde{c})  \right.\nonumber\\
 &\qquad \left. {}
+ A(\tilde{c})x+(1-x)V(\tilde{c})-A(\tilde{b})x-(1-x)V(\tilde{b})|^{p}dx \right)^{1/p}\\
&= \left(\int_{0}^{1}{|A(\tilde{a})x+(1-x)V(\tilde{a})-A(\tilde{b})x-(1-x)V(\tilde{b})|}^{p}dx \right)^{1/p}\\
&= D(\tilde{a},\tilde{b})
\end{align*}
\end{proof}
\end{itemize}

The above arguments proves the metric properties for the proposed distance measure when $V(\tilde{a}) \geq 0$. Similarly it can be proved for the case when $V(\tilde{a})<0$.

\section{Proposed ranking procedure}
\label{sec:4}

\begin{defn}
Let $\mathit{IF}$ denotes the set of all intuitionistic fuzzy numbers. We define a function $\delta : \mathit{IF} \longrightarrow$ $\lbrace-1,1\rbrace $ , $\forall \tilde{a} \in \mathit{IF}$ as 
\begin{equation*}
\delta(\tilde{a})= \left\{
\begin{array}{cc}
1 & , V(\tilde{a}) \geq 0 \\
-1 & , V(\tilde{a}) < 0
\end{array}
\right.
\end{equation*} 
\end{defn}

\begin{defn}
The ranking method based on distance measure is denoted and defined as follows:
\begin{align*}
\rho(\tilde{a})&= \delta(\tilde{a})D(\tilde{a},\tilde{0})\\
&= \left\{
\begin{array}{cc}
d^{(p)}([A(\tilde{a}),V(\tilde{a})],[A(\tilde{0}),V(\tilde{0})]) & ,  V(\tilde{a}) \geq 0 \\
-d^{(p)}([V(\tilde{a}),A(\tilde{a})],[V(\tilde{0}),A(\tilde{0})]) &, V(\tilde{a}) <0
\end{array}
\right.
\end{align*}
 where $\tilde{0}=\langle(0,0,0,0);0,1\rangle$ is considered to be origin.
 
 The steps of distance measure algorithm are:
\begin{itemize}
\item[]Step-I: Compute value index and ambiguity index for each TRIFN ($V(.)$ and $A(.)$)
\item[]Step-II: Compute distance measure between $\tilde{0}$ and each TRIFN ($D(.,\tilde{0})$)
\item[]Step-III: For each arbitrary TRIFNs $\tilde{a}$ and $\tilde{b}$, we define the ranking to be
\begin{enumerate}
\item $\tilde{a} \succ \tilde{b}$ if and only if $\rho(\tilde{a})\succ \rho(\tilde{b})$,
\item $\tilde{a} \prec \tilde{b}$ if and only if $\rho(\tilde{a})\prec \rho(\tilde{b})$ ,
\item $\tilde{a} \sim \tilde{b}$ if and only if $\rho(\tilde{a})\sim \rho(\tilde{b})$,
\end{enumerate} 
\end{itemize}
\end{defn}

\section{Illustrative example}
\label{sec:5}  
\begin{exm} Let us take three sets of intuitionistic fuzzy numbers. 
\\\\Set I: $\tilde{a}=\langle(0.5,0.7,0.9);0.7,0.2\rangle$, $\tilde{b}=\langle(0.2,0.3,0.4);,0.6,0.4\rangle$, $\tilde{c}=\langle(0.4,0.7,0.9);0.6,0.3\rangle$
\\\\Set II: $\tilde{a}=\langle(0.10,0.19,0.25,0.30);0.7,0.2\rangle$, $\tilde{b}=\langle(0.12,0.2,0.23,0.28);0.8,0.1\rangle$, $\tilde{c}=\langle(0.21,0.27,0.32,0.35);0.6,0.3\rangle$
\\\\Set III: $\tilde{a}=\langle(0.2,0.5,0.7);0.7,0.2\rangle$, $\tilde{b}=\langle(0.2,0.3,0.9);0.6,0.4\rangle$, $\tilde{c}=\langle(0.2,0.4,0.5,0.9);0.5,0.3\rangle$
\\\\The table below illustrates the procedure for ranking of intuitionistic fuzzy numbers by proposed method:
\begin{center}
\begin{tabular}{ccccccc}
\multicolumn{7}{c}{} \\
\hline
Sl.no & IFN &  $V(\tilde{a})$ &  $A(\tilde{a})$ & $\rho(.),p=2$ &$\rho(.),p=3$ & Result \\
\hline
\multirow{3}{*}{SetI} &$\tilde{a}$ &0.5250 & 0.0990 &0.3353 &0.3544\\
& $\tilde{b}$ & 0.1800 & 0.0400 & 0.1171 & 0.1232&{$\tilde{b}\prec\tilde{c}\prec\tilde{a}$}\\
&$\tilde{c}$ & 0.4441 & 0.1083 & 0.2927 & 0.3067\\ 
\hline
\multirow{3}{*}{Set II} &$\tilde{a}$ &0.1599 & 0.0799 &0.1221 &0.1241\\
& $\tilde{b}$ & 0.1785 & 0.0623 & 0.1249 & 0.1291&{$\tilde{a}\prec\tilde{b}\prec\tilde{c}$}\\
&$\tilde{c}$ & 0.1885 & 0.0520 & 0.1265 & 0.1319\\ 
\hline
\multirow{3}{*}{Set III} &$\tilde{a}$ &0.3624 & 0.1249 &0.2531 &0.2615\\
& $\tilde{b}$  & 0.2300 & 0.1400 & 0.1868 & 0.1885&{$\tilde{c}\prec\tilde{b}\prec\tilde{a}$}\\
&$\tilde{c}$ & 0.1933 & 0.1800 & 0.1866 & 0.1883\\
\hline
\end{tabular}
\label{Table 1}
%\captionof{table}{Table:1 illustrate the proposed method}
\end{center}
\end{exm}

\begin{exm} Consider the three Intuitionistic fuzzy numbers 
\\$\tilde{a}= \langle (0.3,0.4,0.5,0.6); 0.2,0.4 \rangle$,
\\ $\tilde{b}=\langle (0.1,0.2,0.3,0.4); 0.3,0.5 \rangle $,
\\ $\tilde{c}= \langle (0.5,0.6,0.7,0.8); 0.2,0.6 \rangle$.
\\ By Xu and Yager \cite{xy} which is based on score and accuracy function method the given number are ranked as $\tilde{b}\succ \tilde{a} \succ \tilde{c}$ and apparently the ranking order seems to be $\tilde{c}\succ \tilde{a} \succ \tilde{b}$. But the proposed method indicates that the ranking order among these given set of numbers is $\tilde{c}\succ \tilde{a} \succ \tilde{b}$. Hence, the proposed method is consistent with the human intuition.  
\end{exm}

\begin{exm} Let us take three more Intuitionistic fuzzy numbers 
\\$\tilde{a}=\langle (0.7,0.8,0.9,1.0);0.2,0.5 \rangle$,
 \\$\tilde{b}=\langle (0.3,0.4,0.5,0.6);0.7,0.1 \rangle$,  
 \\$\tilde{c}=\langle (0.5,0.6,0.7,0.8);0.8,0.2 \rangle $.
 \\ According to proposed distance measure, the ranking order of these three numbers are $\tilde{c}\succ \tilde{b} \succ \tilde{a}$. Now negative of these numbers are 
 \\$\tilde{a}=\langle (-1.0,-0.9,-0.8,-0.7);0.2,.05 \rangle$, 
 \\$\tilde{b}=\langle (-0.6,-0.5,-0.4,-0.3);0.7,0.1 \rangle$, and 
 \\$\tilde{c}=\langle (-0.8,-0.7,-0.6,-0.5);0.8,0.2 \rangle $.
\\It is important here to mention that the proposed distance  measure method obeys natural ordering of numbers i.e., for $p=2$, the ranking order of negative of these three numbers are $-\tilde{c}\prec -\tilde{b} \prec \tilde{a}$.
\end{exm}

\begin{exm}
Let us consider a new set of intuitionistic fuzzy numbers with each having same degrees of membership and non-membership as follows:
\\$\tilde{a}=\langle (0.3,0.4,0.5,0.6);0.5,0.3 \rangle$,
 \\$\tilde{b}=\langle (0.7,0.8,0.9,1.0);0.5,0.3 \rangle$,  
 \\$\tilde{c}=\langle (0.2,0.4,0.6,0.8);0.5,0.3 \rangle $.
\\By Novel accuracy function defined by Ye \cite{ye}, the given numbers are ranked as  $\tilde{a}\sim \tilde{b} \sim \tilde{c}$. But according to proposed method, the actual ranking order of these numbers become  $\tilde{b}\succ \tilde{c} \succ \tilde{a}$.
Even with Wei's distance measure method \cite{wei}, we get a similar result.
\end{exm}

\begin{center}
\textbf{Comparison table}
\end{center}
Next we consider the three sets of trapezoidal intuitionistic fuzzy numbers from example 5.2,5.3,5.4. To compare the results obtained by the proposed method with other predefined methods we have presented a comparison table as follows:

\begin{tabular}{lllll}
\multicolumn{5}{c}{}\\
\hline
Authors & Int. fuzzy number & Set I & Set II & Set III\\
\hline
\multirow{3}{*}{Xu and Yager\cite{xy}} & $\tilde{a}$ & 0.6& -0.3& 0.2\\
&$\tilde{b}$ & 0.8& 0.8&0.2\\
&$\tilde{c}$ & 0.8 & 0.6 &0.2\\
\hline
Results& &$\tilde{c}\prec\tilde{a}\prec\tilde{b}$ & $\tilde{a}\prec\tilde{c}\prec\tilde{b}$ & $\tilde{a}\sim\tilde{b}\sim\tilde{c}$\\
\hline
\multirow{3}{*}{Ye\cite{ye}} & $\tilde{a}$ &-0.2& -0.1& 0.3\\
&$\tilde{b}$ & 0.1& 0.5&0.3\\
&$\tilde{c}$ & 0 & 0.8 &0.3\\
\hline
Results& &$\tilde{a}\prec\tilde{c}\prec\tilde{b}$ & $\tilde{a}\prec\tilde{b}\prec\tilde{c}$ & $\tilde{a}\sim\tilde{b}\sim\tilde{c}$\\
\hline
\multirow{3}{*}{Wei\cite{wei}} & $\tilde{a}$ & 0.82&0.70&0.73\\
&$\tilde{b}$ & 0.90 & 0.64& 0.41\\
&$\tilde{c}$ & 0.80 & 0.48& 0.70\\
\hline  
Results& &$\tilde{b}\prec\tilde{a}\prec\tilde{c}$ & $\tilde{a}\prec\tilde{b}\prec\tilde{c}$ & $\tilde{a}\prec\tilde{c}\prec\tilde{b}$\\
\hline
\multirow{3}{*}{D.F.Li\cite{lva}} & $\tilde{a}$ & 0.18&0.29&0.27\\
&$\tilde{b}$ & 0.10 & 0.36& 0.51\\
&$\tilde{c}$ & 0.26 & 0.52& 0.30\\
\hline
Results& &$\tilde{b}\prec\tilde{a}\prec\tilde{c}$ & $\tilde{a}\prec\tilde{b}\prec\tilde{c}$ & $\tilde{a}\prec\tilde{c}\prec\tilde{b}$\\
\hline
\multirow{3}{*}{Li,$\lambda=1/2$ \cite{lra}} & $\tilde{a}$ & 0.16 &0.28&0.26\\
&$\tilde{b}$ & 0.09 & 0.31& 0.48\\
&$\tilde{c}$ & 0.24 & 0.45& 0.25\\
\hline
Results& &$\tilde{b}\prec\tilde{a}\prec\tilde{c}$ & $\tilde{a}\prec\tilde{b}\prec\tilde{c}$ & $\tilde{c}\prec\tilde{a}\prec\tilde{b}$\\
\hline
\multirow{3}{*}{Proposed method,p=2} & $\tilde{a}$ & 0.12&0.18&0.17\\
&$\tilde{b}$ & 0.08 & 0.25& 0.30\\
&$\tilde{c}$ & 0.17 & 0.34& 0.28\\
\hline
Results& &$\tilde{b}\prec\tilde{a}\prec\tilde{c}$ & $\tilde{a}\prec\tilde{b}\prec\tilde{c}$ & $\tilde{a}\prec\tilde{c}\prec\tilde{b}$\\
\hline
\multirow{3}{*}{Proposed method,p=3} & $\tilde{a}$ & 0.04&0.08&0.07\\
&$\tilde{b}$ & 0.02& 0.13& 0.19\\
&$\tilde{c}$ & 0.07 & 0.21& 0.12\\
\hline
Results& &$\tilde{b}\prec\tilde{a}\prec\tilde{c}$ & $\tilde{a}\prec\tilde{b}\prec\tilde{c}$ & $\tilde{a}\prec\tilde{c}\prec\tilde{b}$\\
\hline
\end{tabular}
\label{Table 2}

\section{Conclusion}
\label{sec:6}
In this paper, we have proposed a ranking method based on distance measure and proved that this distance measure satisfies metric properties. Moreover, the advantage of proposed ranking method over other methods is that it can rank intuitionistic fuzzy number with human intuition consistently. We have presented a comparison table to compare the results with other methods.

\section*{Acknowledgment}
The work and research of the first author of this paper is financially supported by TEQIP II,NIT Silchar,Assam,India.


\begin{thebibliography}{222}

\bibitem{al} Allahviranloo,T. Firozja,M.A. Ranking of fuzzy numbers by a new metric. Soft Computing, 14(7), 773-782, 2010.

\bibitem{at} Atanassov, K. Intuitionistic fuzzy sets. Fuzzy Sets and Systems, 20, 87-96, 1986.

\bibitem{atg} Atanassov, K., and Gargov, G. Interval-valued intuitionistic fuzzy sets. Fuzzy Sets and Systems, 31(3), 343-349, 1989. 

\bibitem{chen} Chen, S.J. and Hwang, C.L. Fuzzy Multiple Attribute Decision Making, Springer-Verlag, Berlin, Heidelberg, New York, 1992.

\bibitem{pk} De,P. K., Das.D. Ranking of trapezoidal intuitionistic fuzzy numbers. In Proceedings of IEEE ISDA, 184-188, 2012.

\bibitem{grz} Grzegorzewski,P. Distances between intuitionistic fuzzy sets and/or interval valued fuzzy sets based on the Hausdorff metric. Fuzzy Sets and System,148,319-328, 2004.

\bibitem{hua} Huang,G.S. Liu,Y.S. Wang, X.D. Some new distances between intutionistic fuzzy sets. Proceedings of the Fourth International Conference on Machine Learning and Cybernetics,18-21, 2005. 

\bibitem{hung} Hung, W.L. Yang, M.S. Similarity measures of intuitionistic fuzzy sets based on Hausdorff distance. Pattern Recognition Letters, 25, 1603-1611,2004.

\bibitem{lm} Lakshmana Gomathi Nayagam, V.,Venkateshwari, G. Ranking of intuitionistic fuzzy numbers. In Proceedings of the IEEE international conference on fuzzy systems(IEEE FUZZ2008),1971-1974.

\bibitem{lra} Li,D.F. A ratio ranking method of triangular intutionistic fuzzy numbers and its aaplication to madm problems. Computer and Mathematics with Applications, 60, 1557-1570, 2010.

\bibitem{lva} Li,D.F. Nan,J.X. Zhang,M.J. A Ranking Method of Triangular Intuitionistic Fuzzy Numbers and Application to Decision Making.International Journal of Computational Intelligence Systems, 3(5), 522-530, 2010.

\bibitem{mir} Mitchell, H.B. Ranking intuitionistic fuzzy numbers. International Journal of Uncertainty, Fuzziness and Knowledge Based Systems, 12(3), 377-386, 2004.

\bibitem{kp} Szmidt,E. Kacprzyk,J. On measuring distances between intuitionistic fuzzy sets, Notes IFS 3 , 1–13,1997.

\bibitem{skp} Szmidt,E. Kacprzyk,J. Distances between intuitionistic fuzzy sets, Fuzzy Sets and Systems,114, 505–518,1997.

\bibitem{wan} Wan, S.P., Li,D.F., Possibility mean and variance based method for multi-attribute decision making with triangular intuitionistic fuzzy numbers. Journal of Intelligent and Fuzzy Systems,24(4),743-754, 2013.

\bibitem{wn} Wan, S.P., Power average operators of trapezoidal intuitionistic fuzzy numbers and application to multi-attribute group decision making. Applied Mathematical Modelling, 37, 4112–4126, 2013.

\bibitem{wang} Wang. J. Q., Overview on fuzzy multi-criteria decision making approach. Control Decision, 23, 601-606, 2008.  

\bibitem{wz} Wang, J. Q. and Zhang, Z. Multi-criteria decision making with incomplete certain information based on intuitionistic fuzzy number. Control Decision, 24, 226-230, 2009. 

\bibitem{wli} Wang,L.L., Li,D.F., Zhang,S.S. Mathematical programming methodology for multi-attribute decision making using interval-valued intuitionistic fuzzy sets.Journal of Intelligent and Fuzzy Systems,24(4),755-763, 2013.

\bibitem{wwy} Wang,X.F.,Wang,J.Q.,Yang,W.E. Multi-criteria group decision making method based on intuitionistic linguistic aggregation operators. Journal of Intelligent and Fuzzy Systems,26(1), 115-125, 2014.

\bibitem{wei} Wei, G.W. Some arithmetic aggregation operators with intuitionistic trapezoidal fuzzy numbers and their application to group decision making, Journals of computers,5, 345–351,2010. 

\bibitem{wu} Wu,J.,Cao,Q. Same families of geometric aggregation operators with intuitionistic trapezoidal fuzzy numbers.Applied Mathematical Modelling, 37, 318-327, 2013.

\bibitem{xuc} Xu,Z.,Chen,J. On Geometric Aggregation over Interval-Valued Intuitionistic Fuzzy Information .In Proceedings of the Fourth International Conference on Fuzzy Systems and Knowledge Discovery 2007,2,466-471.

\bibitem{xuzs} Xu,Z.S. Methods for aggregating interval-valued intuitionistic fuzzy information and their application to decision making. Control and Decision, 22(2),  215-219, 2007.

\bibitem{xy} Xu, Z.S., Yager, R.R. Some geometric aggregation operators based on intuitionistic fuzzy sets.International Journal of General Systems, 35,417–433, 2006.

\bibitem{ye} Ye,J. Multicriteria fuzzy decision-making method based on a novel accuracy function under interval-valued intuitionistic fuzzy environment. Expert Systems with Applications, 36, 6899-6902, 2009.

\bibitem{za} Zadeh, L.A. Fuzzy sets. Information and Control, 8(3), 338-356,1965.



\end{thebibliography}
\end{document}